\documentclass[aps,prx,reprint,superscriptaddress]{revtex4-1}
\usepackage[utf8]{inputenc}
\usepackage{amsmath,amsfonts,amssymb}
\usepackage{todonotes}
\usepackage{graphicx}% \graphicspath{{./Figures/}}
\usepackage{braket}
\usepackage[hypertexnames=false,colorlinks=true,linkcolor=blue,citecolor=blue,
            bookmarksopen=false,bookmarks=false,
            pdfstartview=XYZ,pdffitwindow=true,pdfcenterwindow=true]{hyperref}
\usepackage{cleveref}

\renewcommand{\phi}{\varphi}

 \newcommand{\epsi}{\varepsilon}
\newcommand{\eps}{\varepsilon}

\newcommand{\mdrev}[1]{\textcolor{blue}{#1}}

\graphicspath{ {./Figures/} }

\begin{document}

% Use the \preprint command to place your local institutional report number in the
% upper righthand corner of the title page in preprint mode.  Multiple \preprint
% commands are allowed.  Use the 'preprintnumbers' class option to override journal
% defaults to display numbers if necessary
%\preprint{}

\title{Cross-scale excitability in networks of synaptically-coupled quadratic integrate-and-fire neurons}

\author{Daniele Avitabile}
\email[Email: ]{d.avitabile@vu.nl}
\homepage[Home page:]{www.danieleavitabile.com}
\affiliation{Department of Mathematics, Vrije Universiteit Amsterdam,
  De Boelelaan  1111, 1081 HV Amsterdam, The Netherlands}
\affiliation{Inria at Universit{\'e} C{\^o}te d'Azur, MathNeuro Team, 2004 route des Lucioles - BP 93, 06902 Sophia Antipolis Cedex, France}
\affiliation{Amsterdam Neuroscience - System and Network Neuroscience}
\author{Mathieu Desroches}
\email[Email:]{mathieu.desroches@inria.fr}
\homepage[Home page:]{https://www-sop.inria.fr/members/Mathieu.Desroches/}
%\thanks{}
\affiliation{Inria Sophia Antipolis M{\'e}diterran{\'e}e Research Centre, MathNeuro Team, 2004 route des Lucioles - BP 93, 06902 Sophia Antipolis Cedex, France}

\author{G.~Bard Ermentrout}
\email[Email: ]{bard@pitt.edu}
\homepage[Home page: ]{https://www.pitt.edu/~phase/}
\affiliation{Department of Mathematics, University of Pittsburgh, Pittsburgh, PA 15260, USA}

\date{\today}

\begin{abstract}
  From the action potentials of neurons and cardiac cells to the amplification of
  calcium signals in oocytes, excitability is a hallmark of many biological
  signalling processes. In recent years, excitability in single cells has been
  related to multiple-timescale dynamics through \textit{canards}, special solutions
  which determine the effective thresholds of the all-or-none responses.
  However, the emergence of excitability in large populations remains an open
  problem. Here, we show that the mechanism of excitability in large networks and
  mean-field descriptions of coupled quadratic integrate-and-fire (QIF) cells mirrors
  that of the individual components. We initially exploit the Ott-Antonsen ansatz to
  derive low-dimensional dynamics for the coupled network and use it to describe the
  structure of canards via slow periodic forcing.  We demonstrate that the thresholds
  for onset and offset of population firing can be found in the same way as those of
  the single cell. We combine theoretical analysis and numerical computations to develop a novel
  and comprehensive framework for excitability in large populations, applicable not
  only to models amenable to Ott-Antonsen reduction, but also to networks without a
  closed-form mean-field limit, in particular sparse networks.
\end{abstract}
\pacs{}
%\keywords{}

\maketitle

\section{Introduction}\label{sec:intro}
Excitability is a fundamental \textit{all-or-none} property of many living cells
including neurons~\footnote{For simplicity, we describe neuronal excitability,
but this notion extends beyond membrane
biophysics.}. It manifests itself by a very nonlinear response to a sufficiently
strong external input, leading to the emission an action potential before going back
to a rest state, whereas any weaker input has no effect on the cell other than a
small fluctuation of the membrane potential around its equilibrium value. The concept
of excitability is well known to biologists, in particular through the existence of a
non-observable boundary in the response space marking the abrupt transition from rest
to spike. 
However, while the geometry of single-cell excitability is well understood
\cite{rinzel-erm89}, the idea of population excitability (for example in a network of
coupled neurons) has been far less studied. What makes a population respond normally
(such as in working memory tasks \cite{compte00}) or abnormally (such as in seizures
and other pathologies\mdrev{~\cite{kramer2012}}) is a critical question in neuroscience. 

% Most biophysical and phenomenological models of class II membranes \cite{izhikevich}
% ---Hodgkin-Huxley (HH), FitzHugh-Nagumo (FHN), Morris-Lecar (ML)--- share similar
% geometrical features and their excitability threshold is shaped through the intrinsic
% slow-fast dynamics of these models. In particular, in planar models of HH, FHN or ML
% type, and some three-dimensional ones, 
Most class-II membrane models~\cite{izhikevich} ---Hodgkin-Huxley, FitzHugh-Nagumo,
Morris-Lecar--- have a slow-fast structure, with excitability threshold given by
so-called \textit{canard solutions}~\cite{Desroches:2013, Mitry2013,
Wechselberger:2013}.
In single cells, canard solutions underpin complex biological rhythms~\cite{Vo2016},
organise transitions from resting to spiking states~\cite{Moehlis2006}, and from
spiking to bursting regimes~\cite{Kramer2008}.

The canonical class I excitable
systems, such as the QIF neuron models~\cite{ErmentroutKopell:1986}, do not
intrinsically possess multiple timescales.  Nevertheless slow periodic forcing can
bring out a bursting rhythm  in theta neurons and the threshold to bursting dynamics is
again formed by canard solutions~\cite{Desroches:2016}. Networks of QIF neurons
are capable of generating similar bursting rhythms upon periodic
input~\cite{Montbrio:2015,Schmidt2018}, so the question of network excitability in
relation to threshold comes naturally in this context.

In this article, we provide a novel approach to the  question of population
excitability by showing that the geometry of excitability at the microscopic level
scales up to large networks, involving similar key objects related to the slow-fast
nature of the system. Initially, we build on the results by Montbri\'o et
al.~\cite{Montbrio:2015} (extending previous work by Ott and
Antonsen~\cite{OttAntonsen:2008}) in the case of dense (all-to-all) networks of QIF
neurons, with randomly distributed constant inputs following a heavy-tail
(Lorentzian) distribution. For this network, several groups have studied the
existence of a simple mean-field
limit~\cite{Laing2015,Montbrio:2015,devalle2017,DiVolo2018,coombes2019next,Pietras2019}.
We show that excitability in large networks of this type is organised via canards in
the very same
way that it is at the mean-field limit, and we showcase these results
computationally,
by exhibiting an accurate approximation of the network threshold for networks of size
$N\!=\!10^5$. 

What is more, we extend this approach to a much wider class of networks of QIF
neurons, including: networks with heterogeneous weights, sparsely-connected networks,
networks with electrical as well as chemical synaptic coupling, networks with
asymmetrically-spiking neurons, and multi-population networks that combine any of the
above features.

The  geometry of excitability in all such QIF neural networks beautifully
persists across scales, in great generality. For large-enough networks
(and up to the mean-field limit) this persistence reveals itself once we consider the
correct macroscopic variables, namely the firing rate and the mean membrane
potential. This is in contrast to the single-neuron level where the slow-fast
variables endowed with this excitable geometry only encompass the membrane potential.

We show that systems of the type above, at any scale, support a continuous route from
non-bursting to bursting solutions. This continuous route visits canard solutions,
which form an interface for excitable transitions, from
\textit{down network states}  (neural population silent phase) up to \textit{network
bursting} ---as observed in~\cite{Montbrio:2015} but without explanation of the
threshold transition--- as well as for the dual transitions from \textit{up network
states} (neural population tonic firing) down to \textit{network bursting} which was
not been reported before and involves the same canard geometry and dynamics.

In the paper we follow a didactic approach, whereby calculations are performed
initially on the original mean field QIF network derived by Montbri\'o, Paz\'o, and
Roxin in~\cite{Montbrio:2015}, to which we add a synaptic variable as the same group
considered in~\cite{devalle2017}. We henceforth denote this model as the MPR
network. With this example we develop  intuition and all the technical ingredients to
describe the canard population thresholds. We then show how to extend this approach
to more general cases.
The paper is organised as follows: in \cref{sec:MPR} we introduce the MPR model,
present excitability and routes to bursting at single cell and network level; in
\cref{sec:GSPT} we present the mathematical tools to study excitability through
folded-saddle canards, and we use them to interpret network excitability, which we
showcase numerically in networks of $10^5$ neurons; in \cref{sec:GeneralQIF}, we
explain how this approach naturally extends to QIF networks in great generality; 
in \cref{sec:Sparse} we demonstrate that the same continuous routes to bursting
exist in sparse networks, in the absence of an exact mean field limit for the
network; we conclude in \cref{sec:Conclusions}.

\begin{figure*}
  \centering
  \includegraphics{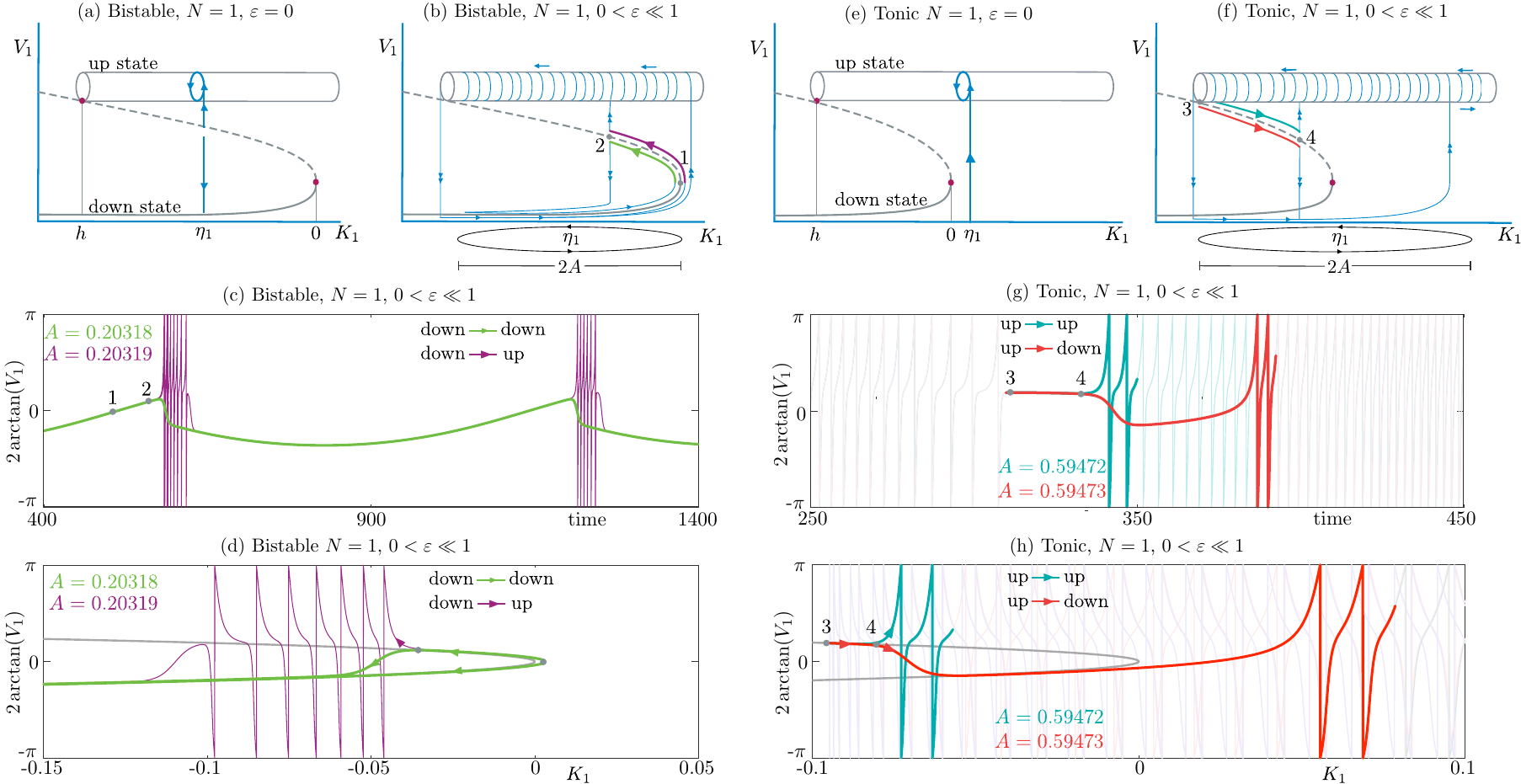}
  \caption{Dynamics of a single QIF neuron ($N=1$ in Eq.~\eqref{eq:newNetw})
  in the bistable regime (a)--(d), and in the tonic regime (e)--(h).
    (a): Sketch (not to scale) of the bifurcation diagram of steady states
    (curve) and periodic solutions
    (cylinder) of the single QIF neuron subject to a
    constant input $K_1 = \eta_1$ ($\eps = 0$). A stable quiescent state (down state)
    coexists with a stable tonic firing solution (up state), separated by an unstable
    equilibrium (dashed curve). A homoclinic bifurcation is present when $K_1 = h$.
    In this intrinsically bistable regime ($K_1 = \eta_1 \in (h,0)$), the cell selects the up or
    down state depending on initial conditions. (b): When $0 <
    \eps \ll 1$, $K_1(t) = \eta_1 + A\sin(\eps t)$ becomes a slowly varying quantity,
    oscillating around the value of $\eta_1$ (ellipses on the $K_1$ axes) with
    amplitude $A$, and
    transitions between the up and the down phases become possible. The
    onset between phases is determined by a family of canard solutions 1--2
    (see text); in the bistable regime they appear in the down-down
    (green), and down-up (purple) transitions. 
    (c): Time profiles of two solutions for the system 
    with slow input $K_1(t)$, displaying a down-down and down-up transition,
    containing a canard segment (1--2). (d) The solutions in (c) are
    plotted in the variables $(V_1,K_1)$, and superimposed on the curve of equilibria
    of the $\eps=0$ system (grey parabola), providing evidence of canard behaviour
    (1--2), and part of the orbits greyed out to enhance visibility. Parameters: 
    $\eps= 0.01$, $J=6$, $\tau_s=0.3$, $\eta_1=-0.2$;
    $A$ values are reported in the panels.
    (e): Sketch of the bifurcation diagram of steady states and periodic solutions
    with constant input ($\eps=0$) in the tonic regime $k_1=\eta_1 >0$). 
    In this regime the cell displays solely the firing solution (up state). (f): When
    $0 < \eps \ll 1$ transitions between the up and the down phases become possible,
    mediated by canard solutions which are possible as up-up and up-down transitions
    (3--4), but not vice-versa. (g): Time profiles of two solutions in the
    tonic regime, with slow input $K_1(t)$, displaying an up-up and up-down transition,
    containing a canard segment (3--4). (h): 
    The solutions in (f) are plotted in the variables $(V_1,K_1)$, and superimposed
    on the curve of equilibria of the $\eps=0$ system (grey parabola), providing
    evidence of canard behaviour (3--4), and part of the orbits greyed out to enhance
    visibility. Parameters: $\eps=0.01$ $J=6$, $\tau_s=0.3$ $\eta_1=0.5$; $A$ values
    are reported in the panels.
}
  \label{fig:sketchN1bi}
\end{figure*}

\section{Population threshold in MPR networks}
\label{sec:MPR}

\subsection{QIF Network model} 
We study a network of $N$ all-to-all coupled QIF neurons.
The $i$th neuron has membrane potential $V_i$, synaptic variable $s_i$, and is
subject to both a background current $\eta_i$, and an external, zero-mean current
$I(t) = A\sin(\epsi t)$, leading to 
%the following evolution equations, 

\begin{equation}\label{eq:newNetw}
  V_i' = V_i^2 + \eta_i +  I(t) + \frac{J}{N} \sum_{j=1}^N s_j,  
\qquad
s_i' =-s_i/\tau_s,
\end{equation}
for $1\leqslant i\leqslant N$.
We refer to the sum $K_i = \eta_i + I(t)$, as the \textit{(external) input} to the
$i$th cell. The ODEs above hold between two consecutive firing times: these are
the finite, computable times at which a membrane potential in the network diverges to
$+\infty$ \footnote{Using the transformation $V_i = \tan \theta_i/2$, one arrives at
$\theta'_i = 1-\cos \theta_i + (1 + \cos \theta_i) \big[ I(t) + \eta_i + \frac{J}{N}
  \sum_{j=1}^N s_j \big],$ $s'_i  =  -s_i/\tau_s$, and computes a firing event when
  $\theta_i$ crosses the value $\pi/ 2$ from below. 
Below we shall present two types of numerical simulations: for $N=1$, we
use this transform, hence thresholds are attained at $V \to \infty$; for $N>1$, we
retain the variables $V_i$, and set a
threshold at $V_i = V_t$, and reset at $V_i=V_r$. The simulation for $N=1$ in
$\theta$ presents no appreciable difference to the one with $N=1$ in $V$.}.
Each time this condition is met by $V_i$ we: (i) stop the simulation, (ii) reset $V_i$
to $-\infty$, (iii) send a spike to \textit{all synapses}, whose values are instantaneously
incremented by an amount $1/N$
;
(iv) restart the simulation of \eqref{eq:newNetw} from
these updated initial conditions. In practice, voltages at $V_t = \infty$, $V_r =
-\infty$ are replaced by finite, large values $V_t = -V_r$.
This system is ideally suited to study excitability
across scales because: (i) we can analyse and compare single cell-dynamics,
$N=1$, network dynamics $N \gg 1$, and mean-field dynamics $N\to \infty$; (ii) one
can switch from constant input currents ($\epsi = 0$) to slowly-varying, oscillatory
currents ($0 < \epsi \ll 1$) in order to uncover transitions between various cellular regimes.

\subsection{Single-neuron excitability} 
Let us set $N=1$, $\epsi=0$, and examine a QIF
neuron with self-coupled synapse~\cite{borgers2018}, subject to a constant
input current $K_1 = \eta_1$, as in~\cref{fig:sketchN1bi}(a,e). When $0 < \tau_s \ll 1$ and
$J\tau_s$ is sufficiently large, the cell supports two \textit{coexisting attracting}
states: an equilibrium (\textit{down state}), and a
periodic solution with tonic
firing (\textit{up state}), separated by an intermediate unstable equilibrium. The
equilibria belong to a curve which folds when the input $K_1$ is null; periodic
solutions collide with unstable equilibria at a homoclinic bifurcation, when 
$K_1=h$.

When $\epsi=0$ and $ \eta_1 \in (h,0)$ in the intrinsically bistable regime,
\cref{fig:sketchN1bi}(a), initial conditions determine whether the voltage is
attracted to the down or to the up state; the threshold is given by the middle
unstable state. When $\eta_1 >0$, the only attractor is the periodic
solution, hence the cell is intrinsically tonic, \cref{fig:sketchN1bi}(e). We are
interested in how the cell (and later on, the network) transitions from the rest
state to the repetitive firing state, and how the two states are concatenated
together to form a bursting state. In a standard QIF model without synapses ($J=0$),
there is no bistability (the up state is to the right of the fold): this changes some
of the waveforms supported by the cell, but not the mechanisms we aim to describe,
namely the bursting transitions between down and up states~\cite{Desroches:2016}.

To study these transitions, one sets $\epsi>0$ small and hence examines slow
forcing~\cite{rinzel-erm89,izhikevich}, as sketched in
\cref{fig:sketchN1bi}(b,f). This causes the input
to oscillate around the mean value $\eta_1$ (see ellipses on
the horizontal axes). The exact dynamics of the system depends on $A$, the amplitude
of the
input oscillations, and on the sign of $\eta_1$. By varying these two parameters, one
can construct a great variety of solutions where up and down states alternate. Some
trajectories stand out, in that they signal the
onset or termination of a phase. With reference to \cref{fig:sketchN1bi}(b),
small-amplitude forcing with average $\eta_1 \in (h,0)$ causes the cell to oscillate around
its rest state (not shown); these are subthreshold oscillations, which stick to the down
state at all times; upon increasing the amplitude $A$ (see ellipse around $\eta_1$),
the trajectories reach a turning point, when $\eta_1+A\!\approx\!0$, that is, for
$A\!\approx\!-\eta_1$; near this value, there are trajectories which follow the
branch of unstable states for increasingly longer times (segments 1--2, termed
\textit{canard segments}), before
jumping to the down state, or to the up state. A temporal profile of solution jumping
down is given in~\cref{fig:sketchN1bi}(c), obtained for $A = 0.20318$, $\epsi = 0.01$
(green); a profile of a solution jumping up is also
in~\cref{fig:sketchN1bi}(c), for $A$ to $0.20319$ (purple curve). The narrow region
of parameter space between $A = 0.20318$ and $A =0.20319$ contains an entire family
of solutions, which surprisingly spend $O(1)$ times near the repelling branch
(segments 1--2).

These counter-intuitive orbits contain canard segments (marked with colors and
numbers in \cref{fig:sketchN1bi}(b)), constitute a computable interface between
subthreshold oscillations (down-down orbits) and bursting states (down-up orbits). 
The canards segments in \cref{fig:sketchN1bi}(b) are marked in the time profiles 
\cref{fig:sketchN1bi}(c), and are also visible in the phase-plane projection 
\cref{fig:sketchN1bi}(d). In the latter the variables $(V_1,K_1)$ are used, and
are superimposed on the curve of equilibria of the $\eps=0$ system (grey parabola),
providing evidence of canard behaviour (1--2).

\begin{figure*} 
  \centering
  \includegraphics[width=\textwidth]{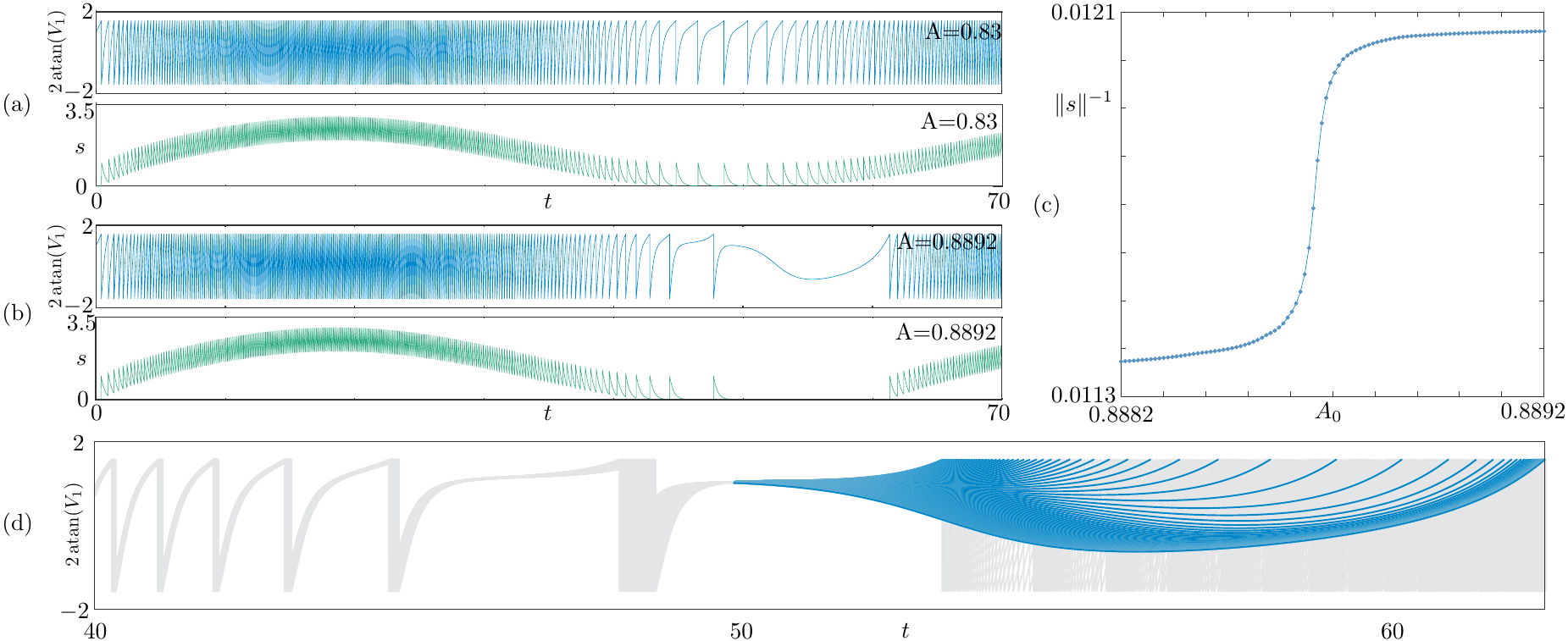}
  \caption{Continuous up-to-down route from non bursting (a) to bursting (b) states
    of an isolated QIF neuron in the tonic regime, upon increasing the amplitude of
    $A$ of the slow forcing. (a): The non-bursting state is visible when the slow
    forcing has small amplitude ($A = 0.83$); the cell exhibits a tonic state with
    slow frequency modulations, because the solution hovers on the top branch of
    \Cref{fig:sketchN1bi}(f), without jumping down. (b): When the forcing amplitude
    is increased slightly ($A=0.8892$) we observe a bursting solution, jumping from
    the top to the bottom branch of \Cref{fig:sketchN1bi}(f). (c): a continuous
    path connects the orbits in (a,b), as $A$ is varied in a narrow band of values;
    the figure shows the reciprocal of the integral $\| s \|$ of $s(t)$ in $t \in
    [0,70]$ as a function of $A$; the sharp increase is typical of canard transitions.
    (d): we plot 100 solutions along the path in (c), near the sharp increase; the
    trajectories (grey) are superimposed, and the canard segments are highlighted in
    blue; with reference to \cref{fig:sketchN1bi}(h), solutions morph from up-up to
    up-down states, while growing a canard segment of type to 3--4.
}

  \label{fig:routeBurstN1}
\end{figure*}

We then move to the intrinsically tonic regime, for $\eta_1 > 0$,
\cref{fig:sketchN1bi}(e--h). We find up-up and up-down orbits with
canards~\cref{fig:sketchN1bi}(f--h), whereas it is possible to prove that down-down
and down-up canards cannot exist, that is, the transition at the fold is a jump. This
is why we have segments of type 3--4, but not of type 1--2, in this scenario.
\cref{fig:sketchN1bi}(g,h) contain many spikes, most of which are greyed-out for
enhancing visibility of the canard segments. In passing we note that segments of type
3--4 are also present in the bistable scenario, but we do not discuss them in the
single cell, for the sake of brevity.

The orbits described above capture excitability transitions at a single-cell level. 
Developing a mathematical understanding of these special orbits is crucial, because
canards act as basin boundaries between different cellular responses: biophysical and
idealised single-cell models support generically \textit{continuous} canard-mediated
transitions \cite{Moehlis2006,Wechselberger:2013}, as we will exemplify in a moment. The transitions are brutally
sharp, but can be continuous, even though they may appear discontinuous upon running
simulations. The main contribution of this paper is to show
that this scenario also occurs generically and robustly in networks of type-I
neurons, of which QIF are universal prototypes.

In addition, the canard-mediated transition from non-bursting to bursting states in
networks of QIF neurons transfers across scales: in an isolated QIF cell, as well as
in a networks of QIF cells, there exists a continuous route from non-bursting to
bursting states, and this path is made of solutions with canard segments such as the
ones seen in \cref{fig:sketchN1bi}. 

\Cref{fig:routeBurstN1} shows an example of such transition in a single cell in the
tonic regime. When $\eps = 0$ the natural solution is a purely tonic one: it locks on
the upper branch of \Cref{fig:sketchN1bi}(e). According to
\Cref{fig:sketchN1bi}(f), we expect to observe a transition from non-bursting to
bursting states, involving solutions with canard segments 3--4, when a slow forcing
$(\eps \neq 0)$ is switched on. When the forcing has small amplitude ($A = 0.83$), the cell exhibits a tonic state
with slow frequency modulations, as seen in \cref{fig:routeBurstN1}(a). This solution
is a non-bursting state, the frequency modulation is present because the solution
hovers on the top branch in \cref{fig:sketchN1bi}(f), which is composed of periodic
orbits with varying period. When  the forcing
is increased slightly, $A=0.8892$, we observe a bursting solution, concatenating a
tonic spiking phase to a quiescent phase; see \cref{fig:routeBurstN1}(b). 

\Cref{fig:routeBurstN1}(c) shows a portion of the continuous path connecting the
non-bursting (\cref{fig:routeBurstN1}(a)) to the bursting
(\cref{fig:routeBurstN1}(b)) state, as $A$ is varied in a narrow band of values. In the
plot, we monitor the states using the reciprocal of the integral of $s$, $\| s \| = \int%_{0}^{70}
s(t) \, dt$ as $A$ varies. The sharp increase along the branch is typical of
canard-mediated transitions, and it is called a \textit{canard explosion}. 
 
In \cref{fig:routeBurstN1}(d) we plot 100 solutions along the path, near the canard
explosion. The trajectories are superimposed in grey, but we highlight the canard
segments in blue: with reference to \cref{fig:sketchN1bi}(h), solutions morph from
up-up to up-down states, while growing a canard segment of type 3--4. 

We label the scenario above as a \textit{up-to-down route to bursting}. Other routes
to bursting are also possible (from down to up states, for instance) when $\eta$
varies. We do not pursue a classification of bursting routes for a single QIF cell,
but we will do so for the mean-field network, after showing that continuous routes to
bursting persist across scales, when $N \to \infty$.

\subsection{Network Excitability} \label{sec:networkExc}
Let us now consider the network \eqref{eq:newNetw}, together with reset
conditions, subject to random background currents: $\eta_i$ are
taken from the Lorentzian distribution with density $g(\eta)=
\Delta/(\pi(\eta-\bar{\eta})^2+\pi\Delta^2)$, hence the network is heterogeneous,
with some neurons in the bistable regime, and others in the tonic regime. However,
the centre of the distribution $\bar{\eta}$ will turn out to play a role. If
$\bar{\eta}<0$ ($\bar{\eta}>0$) we say that distribution, or the network, is
predominantly bistable (tonic). For $N\to\infty$, there is a well-known mean-field
limit~\cite{Montbrio:2015,devalle2017} for the coupled system:
\begin{equation}
\label{eq:oas}
\begin{aligned}
  r' &= \Delta /\pi + 2 rv, \\
  v' &= v^2-\pi r^2 +J s + \bar{\eta} + I(t), \\
  s' &= (-s +r)/\tau_s,
\end{aligned}
\end{equation}
where $r,v,s$ are the mean firing rate, mean membrane potential, and mean synaptic
current, respectively. The external input is now given by $K(t) = \bar \eta + I(t)$.

Note that alternative Ott--Antonsen QIF reductions of QIF networks use amplitude and
phase of a complex-valued order parameter in place of mean voltage and
rate~\cite{Laing2015,coombes2019next}. The order-parameter and rate-voltage
descriptions are related through a conformal mapping~\cite{Montbrio:2015}.

\begin{figure}
  \centering
  \includegraphics{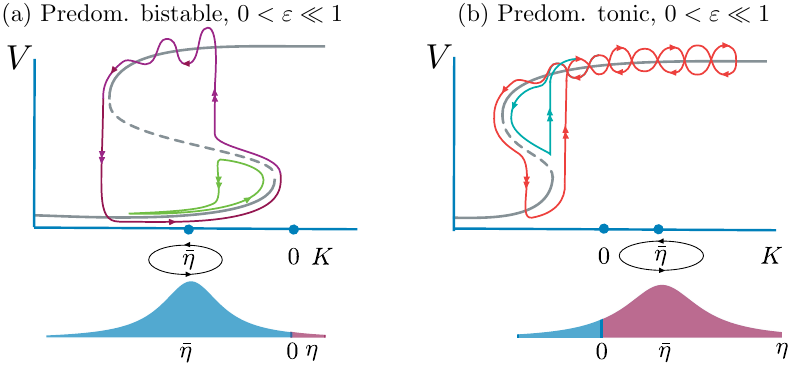}
  \caption{Network with randomly-distributed $\eta_i$ (bottom panels);
    some neurons are bistable (blue), and others tonic (red), with
    distribution centred at $\bar{\eta}$. The $\eps=0$ equilibria lie on $S$-shaped
    curve (grey), whose folds occur for strictly negative values of $K$. The up state
    is now %longer a periodic orbit as in~\cref{fig:sketchN1}, but
    a high-voltage, high-rate equilibrium. The geometry of excitability
    persists when $\eps \ll 1$, and $K(t)$ oscillates slowly around $\bar{\eta}$.
    (a): Down-down and down-up transitions for predominantly bistable distributions
    ($\bar{\eta}<0$). (b): Up-up and up-down transitions for predominantly tonic
    distributions ($\bar{\eta}>0$).
}
\label{fig:sketchInf}
\end{figure}

It is important to remark that the mean-field limit
introduces a new quantity, the population firing rate, $r$~\cite{devalle2017}
that is not a part of the single or finite system of equations. It emerges in the
limit as $N\to\infty$ of the microscopic model \eqref{eq:newNetw}. As seen in
\cref{fig:sketchInf}, when $\eps=0$, hence $K(t) \equiv \bar{\eta}$, the equilibria
of the system lie on an $S$-shaped curve. More precisely, it can be shown that the
curve has no folds or 2 folds \cite[Equation
12 and Figure 2(c)]{schmidt2020bumps}. Henceforth we shall assume that $J$ is
sufficiently large to guarantee the existence of $2$ folds, which must occur for
negative values of $K$.

In the down state all neurons are close to rest (quiescent network state), whereas the up state
corresponds to asynchronous network tonic firing, an averaged version of the tonic
state in~\cref{fig:sketchN1bi}(e--f). The up state can be a stable focus away from
the fold
and have complex eigenvalues.
Between these two stable fixed points is an unstable (saddle) point that serves as a
separatrix between the two stable states. Remarkably, when $0 < \eps \ll 1$, the
geometry of excitability still persists in this macroscopic description, and
transitions are now determined by the distribution peak at $\bar{\eta}$.
We now show that the orbits of the mean field and of the network directly
parallel the transitions of the single neuron model, involving the same canard types.
\begin{figure}
\hspace*{-0.5cm}
  %\centering
  \includegraphics{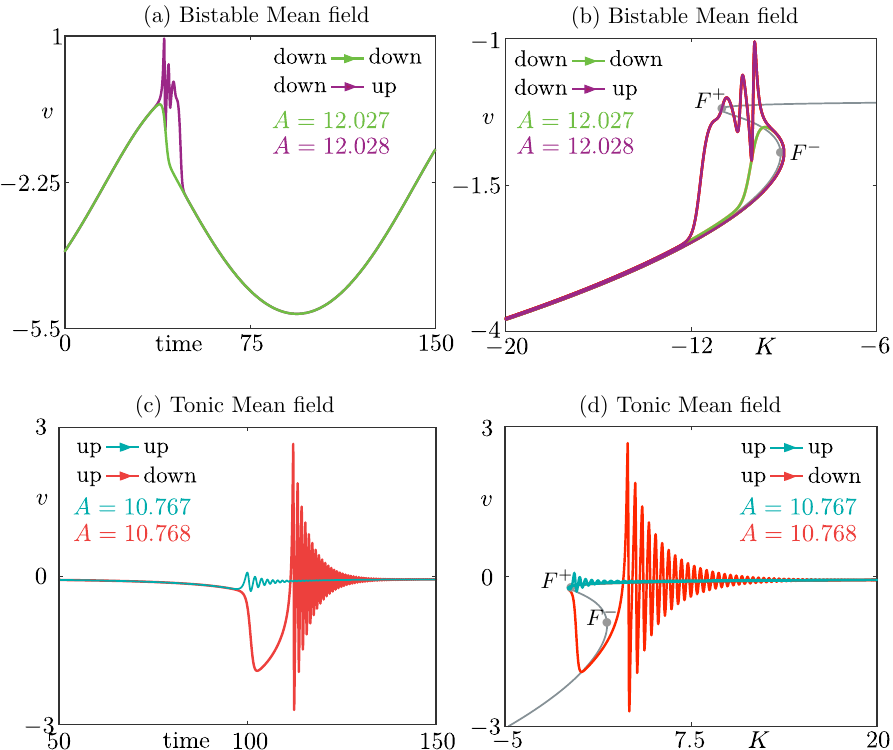}
  \caption{Mean field dynamics of system~\eqref{eq:oas} when the input is slowly
    varying. Top (bottom) panels show the bistable (tonic) regime, reaching
    population bursting via two canard scenarios, down-up and up-down, respectively.
    (a,c): time traces of mean voltage, which display mean-field transitions
    analogous to the single cell ones in \cref{fig:sketchN1bi}(e,g), respectively.
    (b,d), same data as in (a,c), plotted in the variables $(v,K)$, where
    $K=\bar \eta + I(t)$, superimposed on the folded critical manifold (grey); these
    figures are the analogues of \cref{fig:sketchN1bi}(f,h), respectively.
    Parameters are $\Delta=1$, $J=15$, $\tau=0.02$, $\eps=0.05$, and (a,b) $\bar \eta =
    -15.1$, (c,d) $\bar{\eta}=5$. 
  }
  \label{fig:meanFieldFigs}
\end{figure}

\begin{figure*}
  \centering
  \includegraphics{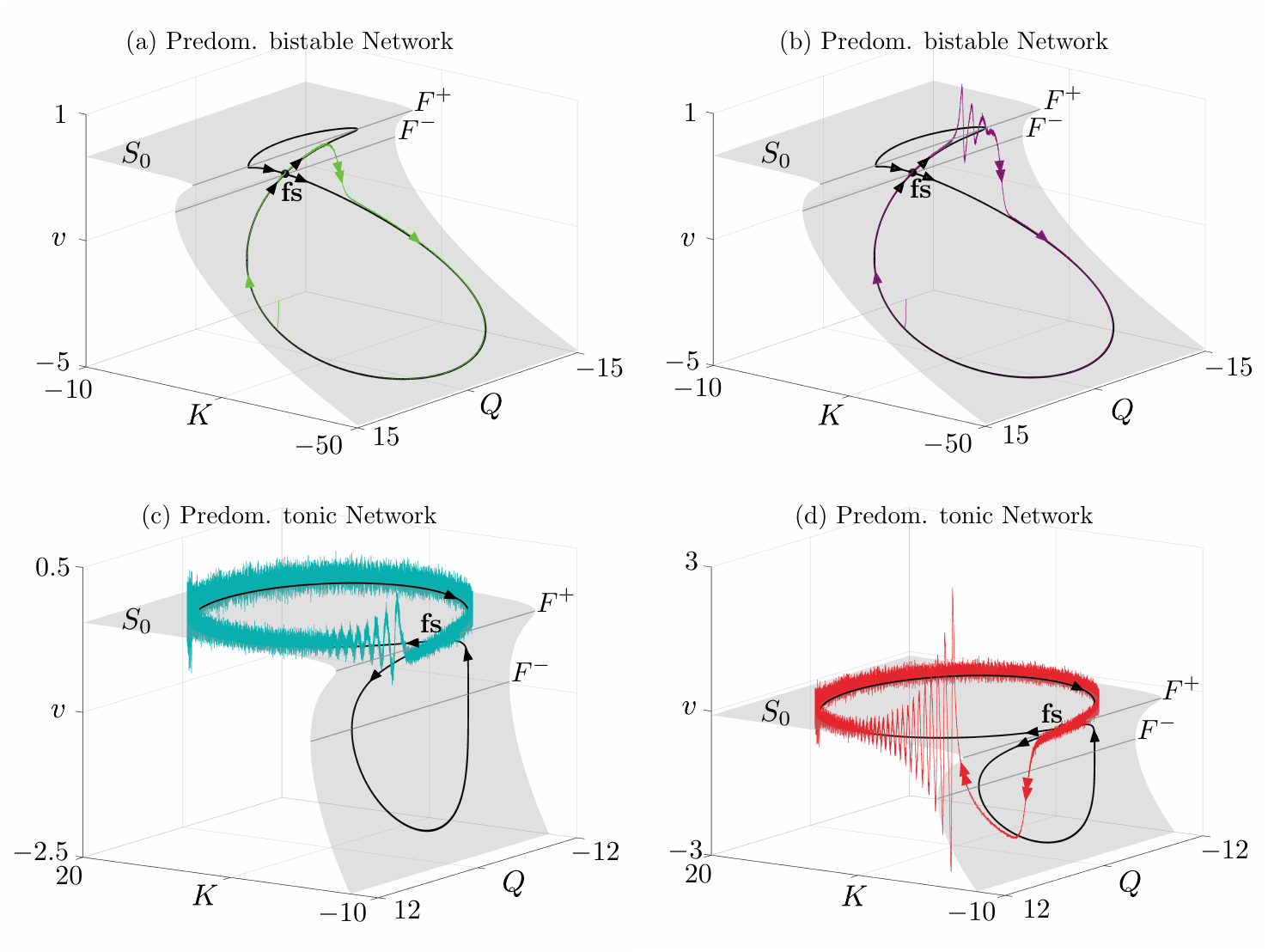}
\caption{Dynamics of network system~\eqref{eq:newNetw} with $N=10^5$ neurons
  when the input is slowly varying, in the bistable (a,b) and tonic (c,d) regimes.
  Shown are the down-down (green) to down-up (purple), and up-up (cyan) to up-down
  (red) transitions in the network, which mirror the mean-field solutions in
  \cref{fig:meanFieldFigs}, as well as the single-cell solutions in
  \cref{fig:sketchN1bi}. To facilitate the comparison between network simulations and
mean-field canard theory (see text), network orbits are shown in the 3D phase space
$(K,Q,v)$, where $Q=I'$. We superimpose them onto the (grey) surface of $\eps =0$
mean-field equilibria, $S_0$, and fold lines $F^{\pm}$ (also shown in
\cref{fig:meanFieldFigs}(b,d)). On the surface are visible the folded saddle
singularity (fs) and its associated canards (black). The networks follows the canard orbits predicted by the mean-field theory remarkably
well.
Parameters are $\Delta=1$, $J=15$, $\tau=0.002$, $\eps=0.05$, $V_t = 100=-V_r$, 
$A$ as in \cref{fig:meanFieldFigs}, and $\bar{\eta}=-15.1$ (a,b),
$\bar{\eta}=5$ (e,f). }
\label{fig:3DFigs}
\end{figure*}

\Cref{fig:meanFieldFigs}(b) shows a simulation with $\epsi =0.05$, $\bar \eta =
-15.1$ for two different values of the maximum amplitude of the stimulus, $A$,
superimposed on the $S$-shaped curve of equilibrium for $\eps = 0$ (grey). One
trajectory (green) follows the down state, along the bottom branch. The trajectory
hugs the unstable branch of fixed points (the canard segment) past the fold $F^-$,
and then jumps down. Small changes in $A$ result in a divergence from the down-down
(green) to the down-up (purple) states. 
% By varying continuously $A$ one can find a
% trajectory that follows the unstable branch up to $F^+$ (not shown).  This trajectory
% acts as an \textit{excitability threshold}, separating down-down from down-up orbits.

The qualitative difference between the two behaviors is more striking in
the time traces of the two curves, shown in \cref{fig:meanFieldFigs}(a),
where the transition to the up state is accompanied by a burst while for the
smaller input, there is just a subthreshold oscillation. The reason for the burst is
that the up state is a stable spiral for a range of input values. 
% We will later show
% that this transition is determined by specific type of canard, the so-called
% \textit{folded-saddle canard}, just as it was in the single neuron case.       
%
When $A \gg -\bar{\eta}$, as expected, the network jumps (without canards) from
asynchronous to synchronous firing, following the upper branch of the $S$-shaped
curve in \cref{fig:sketchInf}(a), and \cref{fig:meanFieldFigs}(b). 

In conclusion, the trajectories displayed in \cref{fig:meanFieldFigs}(a,b) are the
mean-field equivalents of the single cell ones in the bistable regime, shown in
\cref{fig:sketchN1bi}(c,d), respectively. They undergo similar transitions from
down-down to down-up states, with canards being the threshold.

The tonic case ($\bar{\eta} = 5$) is simulated in \cref{fig:meanFieldFigs}(c,d).
A similar situation occurs in this case, but the dynamics revolves around the upper
fold $F^+$, and features canards in up-up and up-down orbits.

\section{Folded-saddle canard behavior across scales} 
\label{sec:GSPT}
We now derive three central results of the paper: firstly, we characterise the
mean-field transitions
described above, valid at the ODE level (system \eqref{eq:oas}), using standard methods
from Geometric Singular Perturbation Theory~\cite{J95}; secondly, we use this
characterisation to infer the existence of canard behavior at the network level
(system \eqref{eq:newNetw}) with $N=10^5$ neurons, and provide numerical evidence of
this phenomenon; thirdly, we explore canard-mediated routes to bursting at the
network level, a feature that persists from the single neuron level.
The latter results are remarkable and novel, as canard behavior in large networks is
greatly unexplored, in particular for systems with resets and random data.

As anticipated, we initially present the theory for the network described above, and
then adapt these results to more general networks.
To study the behavior of~\cref{eq:oas} for small $\epsi>0$, we extend the system
with two ODEs describing the oscillatory dynamics of $K(t)$, namely 
\begin{equation}\label{eq:forcing}
 \begin{aligned}
 K' &=  ~~\eps Q,\\ 
 Q' &= -\eps (K - \bar\eta).
 \end{aligned}
\end{equation}
Note that in~\cite{Montbrio:2015,Pietras2019}, the mean-field limit model
assumes instantaneous synaptic processing, which amounts to taking $\tau_s=0$ and
replacing $s$ by $r$ in the equations for $r$ and $v$. However, assuming
$0<\tau_s\ll1$, that is, a fast synapse, does not change anything about the threshold
analysis below while making it more general.

Now that we have expressed the slow dynamics of the current $I$ using a second-order harmonic
equation, we rescale time in~\cref{eq:oas,eq:forcing} so as to parametrise them by the slow time $\tau = t/\epsi$, as was done in~\cite{Avitabile:2017}, and obtain
\begin{equation}\label{eq:meanFieldExt}
  \begin{aligned}
   \epsi\dot r &= \Delta /\pi + 2 rv, \\
   \epsi\dot v &= v^2-\pi r^2 +J s + K, \\
   \epsi\dot s &= (-s +r)/\tau_s, \\
          \dot K & = Q, \\
	  \dot Q & = -(K-\bar{\eta}).
  \end{aligned}
\end{equation}

To shed further light onto the transition from low-rate (down) states to high-rate
(up) states of the mean-field, it is key to consider the slow limit
of~\cref{eq:meanFieldExt} as the forcing speed $\epsi$ tends to $0$.
Therefore, we set $\epsi=0$ in~\cref{eq:meanFieldExt} and obtain the three algebraic
constraints 
\[
s = r, \qquad r = - \frac{\Delta}{2\pi v}, \qquad  v^2-\pi r^2 +J s + K = 0,
\]
In the $\epsi \to 0$ limit, the system variables $(r,v,s,K,Q)$ evolve on a
three-dimensional manifold in $\mathbb{R}^5$, the so-called \textit{critical
manifold}, given by
\[
  S_0 = \left\{s=r=-\frac{\Delta}{2\pi v}, \quad 0 = K + \psi(v) \right\},
\]
where
\begin{equation}\label{eq:PsiDef}
  \psi(v)=v^2 - \Delta^2/(4\pi v^2) -J\Delta/(2\pi v).
\end{equation}
The subscript $0$ in $S_0$ refers to the fact that this manifold is found
by setting $\epsi=0$ in \eqref{eq:meanFieldExt}. The transitions discussed in this
paper occur when $S_0$ is folded, and it is the union of two attracting and one
repelling submanifolds. These conditions occur
generically: for common choices of parameters, one finds that $S_0$ has two loci
of folds (two lines of folds), $F^+$ and $F^-$, corresponding to the set
$\{D\psi(v):=\psi'(v)=0\}$. A projection of the manifold $S_0$ onto the $(K,v)$ plane is visible
in~\cref{fig:meanFieldFigs} (b,d), where the fold lines project onto points $F^{\pm}$; compare with~\cref{fig:3DFigs} where $S_0$ is projected onto the $(K,Q,v)$ space and the fold lines are fully visible.

The $\epsi\to0$ limit introduced above corresponds to a differential-algebraic
problem referred to as the \textit{slow subsystem} of the original equation, and in
the present case it reduces to
\begin{equation}\label{eq:slowsub}
 \begin{aligned}
    0 &= K+\psi(v), \\
    \dot K &= Q, \\
    \dot Q &= -(K-\bar\eta).
 \end{aligned}
\end{equation}

The algebraic constraint in~\cref{eq:slowsub} hides the dynamics of $v$ in this slow
limit. To reveal it, we then differentiate the constraint with respect to time, and
obtain the following set of ODEs defined for $(v,K,Q)\in S_0$
\begin{equation}\label{eq:slowSubs}
  \begin{aligned}
        -\psi'(v)\dot{v} & = Q, \\
  \dot Q             & = \bar{\eta}+\psi(v).
  \end{aligned}
\end{equation}

One can relate \cref{eq:slowSubs} to the canards shown in \cref{fig:meanFieldFigs}:
we are considering the $\epsi=0$ dynamics, hence we focus on the black curves on
$S_0$, and we take the one in \cref{fig:meanFieldFigs}(a) as an example. It would
appear that \cref{eq:slowSubs} breaks down at points where $\psi'(v)=0$, that is,
along the fold set $F^\pm$. Along such folds the first equation reduces to $Q=0$
hence \cref{eq:slowSubs} is undefined at points along the folds where $Q\neq
0$. However, an inspection of
\cref{fig:meanFieldFigs}(a) shows that the flow is well defined at one specific
point, which is called a folded singularity and is marked as
(fs)~\footnote{In fact, there are two other points on $F^+$ where the trajectory seems
  to cross the fold. However, these are not associated with canard dynamics, and the
flow is not passing through the fold, as the system \cref{eq:slowSubs} is singular at
those points.}. 

It is around this point that canard solutions are born, because the trajectory
passes through (fs) from an attracting to a repelling sheet of $S_0$. We also
note that for $\epsi \ll 1$ this behaviour persists (green curve in
\cref{fig:meanFieldFigs}(a)): slow-fast theory predicts~\cite{krupa2001}
that canards at $\eps=0$ survive for small-enough $\eps>0$ and, as we will see below,
they organise the excitable structure of QIF networks.

We therefore investigate the passage through (fs) more precisely: intuitively, at
(fs) $\psi'(v)=0$ and $Q=0$, so that the quotient $Q / \psi'(v)$ stays finite, and
the flow of \cref{eq:slowSubs} is well defined. We formalise this step by: (i)
desingularising \cref{eq:slowSubs} with a time rescaling, (ii) identifying (fs) as an
equilibrium point of the desingularised problem, (iii) classifying the type of
equilibrium, which in turn determines the type of canard in the original system (as
done in \cite{Avitabile:2017}).

In step (i), we desingularise \cref{eq:slowSubs}, and rescale time by a factor
$-\psi'(v)$, which eliminates the prefactor to $\dot v$, and regularises the problem,
leading to the \textit{desingularised reduced system (DRS)}\footnote{An important
  subtlety is that the time rescaling by $-\psi'(v(t))$ depends on the state variables
  $v$, hence the time orientation depends on the position on $S_0$. In fact, the time
  rescaling transforming the reduced system \cref{eq:slowSubs} into the desingularised reduced system
\cref{eq:DRS} is such that the flow in both systems has the same orientation on
the attracting sheets of $S_0$ but opposite orientation on its repelling sheet.}
\begin{equation}\label{eq:DRS}
  \begin{aligned}
        v' & = Q, \\
       Q' & = -\psi'(v)\left(\bar{\eta}+\psi(v)\right).
  \end{aligned}
\end{equation}

In step (ii) we look for equilibria of \cref{eq:DRS}, which satisfy $\psi'(v) = 0$
and $Q=0$, hence they geometrically coincide with (fs). Such equilibria are of the
form $(v,Q) = (v_*,0)$, where $v_*$ satisfies $\psi'(v_*) = 0$.

In step (iii) we study the linear stability of these
equilibria, which is determined by the Jacobian matrix 
\[
  J=\begin{bmatrix}
    0                      &  1 \\
    -\psi''(v_*) \psi(v_*) & 0
  \end{bmatrix},
\]
whose eigenvalues are given by
\[
  \lambda = \pm \sqrt{-\psi''(v_*) \psi(v_*)}.
\]
The equilibrium $(v_*,0)$ is therefore either a saddle or a center: in the former
case the point (fs) is called a \textit{folded saddle}, and gives rise to folded-saddle canards; in the latter case (fs) is a \textit{folded centre}; it is
known that this singularity does not give rise to canards, but rather to a
discontinuous transition. 
\begin{figure*}
  \centering
  \includegraphics{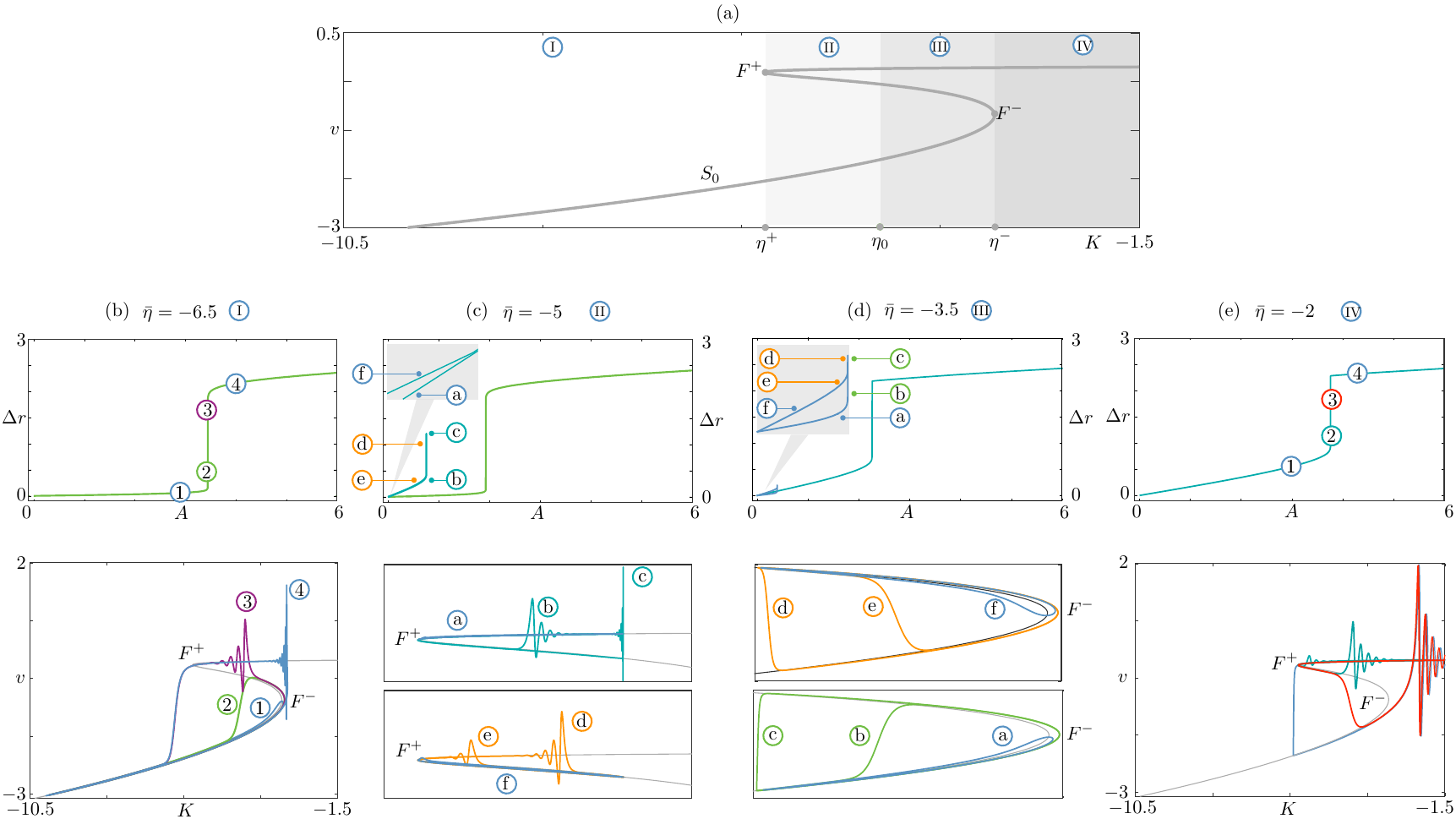}
  \caption{Continuous routes to bursting and non-bursting solution branches, in the mean field model~\eqref{eq:meanFieldExt} for four different values of $\bar{\eta}$ corresponding to the four different scenarios described in the main text. The associated four regions I, II, III and IV are highlighted on top of the critical manifold $S_0$ in panel (a). The three values of $\bar{\eta}$ at the transitions between these scenarios are denoted $\eta_+:=\eta(F^+)$, $\eta_-:=\eta(F^-)$ and $\eta_0=(\eta^++\eta^-)/2$, where $F^\pm$ are the two folds of $S_0$. The chosen values of $\bar{\eta}$ in each region are $\bar{\eta}=-6.5$, $\bar{\eta}=-5$, $\bar{\eta}=-3.5$ and $\bar{\eta}=-2$, respectively. Panels (b--d) display both the solution branches obtained by varying the forcing amplitude $A$ for a given value of $\bar{\eta}$ (top), and a selection of solutions on the branch, plotted in the phase plane $(K,v)$ on top of $S_0$ (bottom). As observed, a continuous branch of solutions bridging from the non-bursting regime to the bursting regime always exists, regardless of the $\bar{\eta}$ value. This branch connects in parameter space a down-down non-bursting solution to a down-up bursting one for $\bar{\eta} < \eta_0$ (panels b,c), or an up-up non-bursting solution to an up-down bursting one for $\bar{\eta} > \eta_0$ (panels d,e). Additionally, for $\eta^+<\bar{\eta}<\eta^-$, another solution branch exists, starting at low $A$ amplitude and which does not connect to the bursting regime; this branch contains up-up (resp. down-down) solutions for $\eta^+<\bar{\eta}<\eta_0$ (resp. $\eta_0<\bar{\eta}<\eta^-$). Parameter values are: $\Delta = 1$, $J=15$, $\tau_s=0.02$, $\epsi=0.05$; $\bar{\eta}$ and $A$ as indicated in the panels.}
  \label{fig:scenarios}
\end{figure*}

A quick calculation reveals that: (1) In the bistable regime, when $\bar{\eta}<0$,
there are two (fs) points, one on $F^+$ and one on $F^-$, both of which are folded
saddles; both folded saddles are visible (in projection) as $F^+$, $F^-$ in
\cref{fig:meanFieldFigs}(b), whereas only the latter is shown in
\cref{fig:3DFigs}(a,b). (2) In the tonic regime, when $\bar{\eta} > 0$, the
singularity on $F^+$ is still a folded saddle, whereas the one on $F^-$ is a folded
centre, and both are visible in projection in \cref{fig:meanFieldFigs}(d); we show only the folded
saddle as (fs) in \cref{fig:3DFigs}(c,d). 
In passing, we note that these findings remain valid for a much larger class of networks and
mean field limits, as we will show below.

\Cref{fig:3DFigs} displays the dynamics of~\cref{eq:slowSubs} ($\epsi = 0$, in black) superimposed on
the dynamics of~\cref{eq:newNetw} with $N=10^5$ ($\epsi >0$, in color), for both
subthreshold and suprathreshold forcing. \Cref{fig:3DFigs}(a) shows, in
green, the \emph{full network orbit with $N=10^5$} for $A$ slightly below
threshold. A canard segment is visible, where the trajectory hugs the black curve on
the fold before falling back to the down state. \Cref{fig:3DFigs}(b) shows the
same projection, but for slightly larger $A$; in this case, the trajectory makes
the jump to the up state before falling back down. The folded-saddle canards
predicted by mean-field theory (in black) are in striking agreement with the
behaviour of the network, and the correspondence with the similarly coloured
mean-field transitions in~\cref{fig:meanFieldFigs}(a,b) and
\cref{fig:sketchN1bi}(c,d) is remarkable.
These solutions define analogue scenarios to the down-down and down-up states for the
single neuron model, but they are exhibited at the level of the network.
Interestingly, at single cell level we have a bistable neuron ($\eta_1 <0$), while at
network level some neurons will be tonic, but the network is predominantly bistable
($\bar{\eta} < 0$).
Similar considerations apply to the predominantly tonic network scenario, showcased in
\cref{fig:3DFigs}(c,d), and corresponding to \cref{fig:meanFieldFigs}(c,d) and
\cref{fig:sketchN1bi}(g,h). 

\subsection{Continuous routes to bursting}
We will now show that for any value of
$\bar{\eta}$, the network robustly supports a continuous route from a non-bursting
state to a
bursting state upon increasing the input amplitude; this transition involves a canard
explosion, as seen in \cref{fig:routeBurstN1} for the single cell.

For sufficiently large coupling values $J$, the critical manifold $S_0$ is
S-shaped, with folds occurring at $\eta_+ < \eta_- < 0$. 
As depicted in \cref{fig:scenarios}(a), there exist $4$ different scenarios depending
solely on the value of $\bar{\eta}$ with respect to the fold values $\eta_\pm$
and the midpoint $\eta_0 = (\eta_+ + \eta_-)/ 2$ between them.

\textit{Case I: $\bar{\eta} < \eta_+$.} Since the forcing oscillates around
$\bar{\eta} < \eta_+$, when a low-amplitude forcing is switched on, the
network can only oscillate near the bottom branch of $S_0$ (orbit 1 in
\cref{fig:scenarios}(b)). Upon increasing the amplitude of the forcing, the orbit
reaches the fold and hugs the middle unstable branch of $S_0$ (orbit 2 in
\cref{fig:scenarios}(b)). A continuum of orbits with canard segments are now visited
by the network: in this transition to bursting visible in the $(A,\Delta r)$ diagram
in \cref{fig:scenarios}(b), the branch of solutions is almost vertical, as $A$ varies
in a tiny region of parameter space. As we ascend the branch, we pass from a
down-to-down solution (green, labelled as $2$) to a down-to-up solution (purple,
labelled $3$) with canard segments. Past the vertical branch, we obtain a fully
developed bursting solution ($4$). The solution branch in \cref{fig:scenarios}(b) 
reveals a \textit{continuous down-to-up route from non bursting to bursting network
solutions}. All solutions with canard segments in this case are of down-to-up or
down-to-down type, as witnessed by the green and purple coloring.

\textit{Case II: $\eta_+ < \bar{\eta}< \eta_0$.} Since we enter the bistable region of
$S_0$, we can switch on the forcing near two starting points, on the lower and upper
branch, respectively. If one starts from the lower branch, the same considerations
as in Case I are valid, and we have a continuous down-to-up route (see green branch
in \cref{fig:scenarios}(c)). Starting from the upper branch, low amplitude forcing
generate state hovering near the upper branch (solution a in
\cref{fig:scenarios}(c))). When the forcing amplitude increases the orbit grows a
canard segment from the upper fold $F_+$ (solution b). However, due to the proximity
between $\bar{ \eta}$ and $F^+$, the branch can only grow canards up to
the solution labelled c, and folds back onto itself while displaying solutions d--f.
Therefore in Case II there is a continuous down-to-up transition (and no up-to-down)
transition to bursting.

\textit{Case III: $\eta_0 < \bar{\eta} < \eta_-$.} This scenario is the mirror image
of Case II. The network possesses \textit{continuous down-to-up transition to
bursting}, but no down-to-up transition, which is interrupted
(\cref{fig:scenarios}(d)). The continuous transition is explained in case IV below. 

\textit{Case IV: $\bar{\eta} > \eta_-$.} When the forcing is small, the solution
can only stay near the upper branch of $S_0$ (solution 1 in \cref{fig:scenarios}(e)).
We can still transition continuously from this non-bursting solution to a bursting
solution (solution 4) via canards that start near $F_+$. This case mirrors Case I,
but involves canards of up-to-up and up-to-down type.

% \textit{Limiting cases.} We have numerically explored the cases $\bar{\eta }=
% \eta^\pm$, in which observe that the respective non-bursting branch shrinks and
% disappears, while the bursting one persists. Also, in the case
% $\bar{\eta}=\eta_0$, the down-to-up continuous route persists, while the
% up-to-down branch is disconnected.

\section{Extension to general QIF networks}
\label{sec:GeneralQIF}
\subsection{Networks with heterogeneous currents}
Let us now consider generalisations of the MPR network~\cref{eq:newNetw} and
mean-field limit \cref{eq:oas}, for which the excitability scenario described above
still holds. We will discuss the generalisations only at the level of the mean field,
and refer to existing literature for descriptions of the corresponding microscopic
networks.

The starting point is the following generalisation of the MPR mean field~\eqref{eq:oas}
\begin{equation}
 \label{eq:oasgen}
  \begin{split}
    r' & = \Delta/\pi + 2rv + (\Gamma/\pi-g)s \\
    v' & = v^2 -\pi^2 r^2 + (J+g \ln a) s + \bar{\eta} + I(t) \\
    \tau_s s' & = (-s+r)
  \end{split}
\end{equation}
where: $I(t) = A\sin(\eps t)$ is a slow zero-mean periodic forcing as before;
$\Delta$, $J$, $\bar{\eta}$ and $\tau_s$ are parameters as before; $\Gamma$, $g$, and
$a$ are additional parameters. This generalisation encompasses a variety of exact
mean-field limits of QIF networks including:
\begin{enumerate}
  \item The MPR network with
    heterogeneous background currents $\eta_i$ sampled using a Cauchy distribution
    with peak at $\bar{\eta}$ and half-width at half-maximum (HWHM) $\Delta$. To
    recover this model from~\cref{eq:oasgen}~\cite{Montbrio:2015,DiVolo2018} one sets
    $\Gamma=0$, $g = 0$, $\tau_s = 0$.
  \item The MPR network with heterogenous (all-to-all) synaptic
    coupling~\cite{Montbrio:2015}, obtained for $\Gamma\neq0$, $g=0$, $\tau_s=0$.
  \item The MPR network with first-order fast or slow
    synapses~\cite{Schmidt2018,devalle2017} characteristic time $\tau_s$, which
    corresponds to $\Gamma=0$, $g=0$, $\tau_s\neq0$. Note that the analysis done
    above on~\cref{eq:newNetw} assumed $0<\tau_s\ll1$ (fast synapse), however it is
    still valid for $\tau_s=O(1)$ (slow synapse), because the folded-saddle structure
    only requires the existence of a slow periodic forcing~\eqref{eq:forcing}. This dynamics
    also persists with second-order synapses (data not shown).
  \item The modified QIF network from~\cite{Pietras2019} with electrical coupling,
    which corresponds to $\Gamma=0$, $g\neq0$, $a=1$, $\tau_s=0$.
  \item The modified QIF network studied in~\cite{montbrio2020} with electrical
    coupling and asymmetric spikes, which differs from the previous case only by
    $a\neq1$.
\end{enumerate}

Before further extending the networks analysable with the proposed formalism, let us
rewrite~\cref{eq:oasgen} using the generalised coefficients $\tilde \Gamma (\Gamma,g) =
\Gamma / \pi -g$, $\tilde J(J,g,a) = J + g \ln a$, yielding
%\[
%  \begin{aligned}
%    r' & = \Delta/\pi + 2rv + \bar \Gamma s       && := \rho(r,v,s) \\
%    v' & = v^2 -\pi^2 r^2 +  \eta + \bar J s + I(t) && := \nu(r,v) + \bar J s + I(t) \\
%    s' & = (-s+r)/\tau_s                          && := \sigma(r,s),
%  \end{aligned}
%\]
\begin{equation}\label{eq:meanFieldGen}
  \begin{aligned}
    r' & = \Delta/\pi + 2rv + \tilde \Gamma s := \rho(r,v,s) \\
    v' & = v^2 -\pi^2 r^2 +  \bar{\eta} + \tilde J s + I(t) := \nu(r,v) + \tilde J s + I(t) \\
    s' & = (-s+r)/\tau_s := \sigma(r,s).
  \end{aligned}
 \end{equation}
It is apparent that one can analyse the slow-fast structure of this system in the
exact same way as we have done for~\cref{eq:meanFieldExt}, with a cubic-shaped
critical manifold given by 
$$S_0 = \left\{s=r=-\Delta/(\pi(2v+\tilde \Gamma)),\;0 = K + \psi(v)\right\},$$
where $\psi$ is as in \cref{eq:PsiDef} with $\tilde J$ instead of $J$. It is apparent
that the same folded-saddle dynamics organise the excitable structure of the
corresponding mean-field model and can be observed (data not shown) in associated
large-enough generalised QIF networks. 

Instead of pursuing this analysis, we introduce first a further generalisation,
namely we consider the exact mean-field limit of $p$ synaptically coupled
populations of QIF networks, where we suppose, for simplicity, that only one
population (the $k$-th one) receives a slow external periodic forcing. The coupled equations read
\begin{equation}\label{eq:meanFieldExtppop}
  \begin{aligned}
   \epsi\dot r_i &= \rho_i(r_i,v_i,s_i), \\
   \epsi\dot v_i &= \nu_i(r_i,v_i)+\sum_{j=1}^p\tilde J_{ij}s_j+K\delta_{ik}, \\
   \epsi\dot s_i &= \sigma_i(s_i,r_i), \\
          \dot K & = Q, \\
	  \dot Q & = -(K-\bar{\eta}_k),
  \end{aligned}
\end{equation}
for  $i=1,\cdots,p$, where $\delta$ is the Kronecker symbol, and $\rho_i$, $\nu_i$,
$\sigma_i$ are the functions defined in~\cref{eq:meanFieldGen} for
population-specific choices of parameters $\tilde \Gamma_i$, $\Delta_i$,
$\bar{\eta}_i$, and $(\tau_s)_i$
. The critical
manifold of~\cref{eq:meanFieldExtppop} is defined by the algebraic constraints
\begin{multline}
\Psi_i(v_1,\cdots,v_p)
:=\nu_i\left(-\frac{\Delta_i/\pi}{2v_i+\tilde{\Gamma}_i},v_i\right) \\
-\sum_{j=1}^p\tilde{J}_{ij}\frac{\Delta_j}{\pi(2v_j+\tilde{\Gamma}_j)}=0.
\end{multline}
where
\begin{equation*}
r_i =s_i = -\frac{\Delta_i/\pi}{2v_i+ \tilde \Gamma_i}
\end{equation*}
for $i=1,\ldots,p$. Hence, the critical manifold $S_0$ can be compactly written as
\begin{multline}
S_0=\bigg\{(v_1,\cdots,v_p,K,Q)\in\mathbb{R}^{p+2} \colon \\ 
0=\Psi_i(v_1,\cdots,v_p)+K\delta_{ik}, \quad i = 1,\ldots,p\bigg\}.
\end{multline}
The expression for $S_0$ contains $p$ independent algebraic conditions in
$\mathbb{R}^{p+2}$, therefore the critical manifold is indeed a surface, which is
consistent with the fact that~\cref{eq:meanFieldExtppop} has two slow variables. 
As a consequence, one can write the reduced system associated
with~\cref{eq:meanFieldExtppop} in the form
\begin{equation}\label{eq:RSmeanFieldExtppop}
 \begin{aligned}
 0 &= \Psi_i(v_1,\cdots,v_p)+K\delta_{ik}, \qquad i=1,\ldots,p\\
 \dot{K} &= Q\\
 \dot{Q} &=-(K-\bar{\eta}_k)
 \end{aligned}
\end{equation}
The system above mirrors the constrained system~\eqref{eq:slowsub} in the
single-population MPR network. Proceeding like in the single-population case, we
differentiate the algebraic constraints with respect to time and project the
resulting limiting system onto the $(v_k,Q)$-plane to obtain
\begin{equation}\label{eq:RSbismeanFieldExtppop}
 \begin{aligned}
 -\partial_{v_k} \Psi_k(v_1,\cdots,v_p)\dot{v}_k &=Q,\\
  \dot{Q} &=\bar{\eta}_k+\Psi_k(v_1,\cdots,v_p),
 \end{aligned}
\end{equation}
where $v_1,\cdots,v_p$ satisfy the algebraic constraints defining the critical
manifold, that is, the first $p$ equations in~\eqref{eq:RSmeanFieldExtppop}.
We recognise in~\cref{eq:RSbismeanFieldExtppop} the same form as \cref{eq:slowSubs}
for the one-population mean-field limit. We note that the starting system has $p+2$
equations, but the reduced system \eqref{eq:RSbismeanFieldExtppop} has only $2$
equations, and it is singular at the fold set of the $k$th population system, given by
the condition $\partial_{v_k} \Psi_k(v_1, \cdots, v_{k*}, \cdots,v_p) =0$. 

The 2-dimensional system \cref{eq:RSbismeanFieldExtppop} is singular along that fold,
and it can be desingularised as in the one-population problem. The folded-saddle and
folded-centre classification carries through in this case. Hence we can conclude that
the same canard-induced excitability scenario appears in the generic $p$-population
case described above.

\section{Canard transitions across scales in sparse networks}
\label{sec:Sparse}
\begin{figure*}
  \centering
  \includegraphics[width=\textwidth]{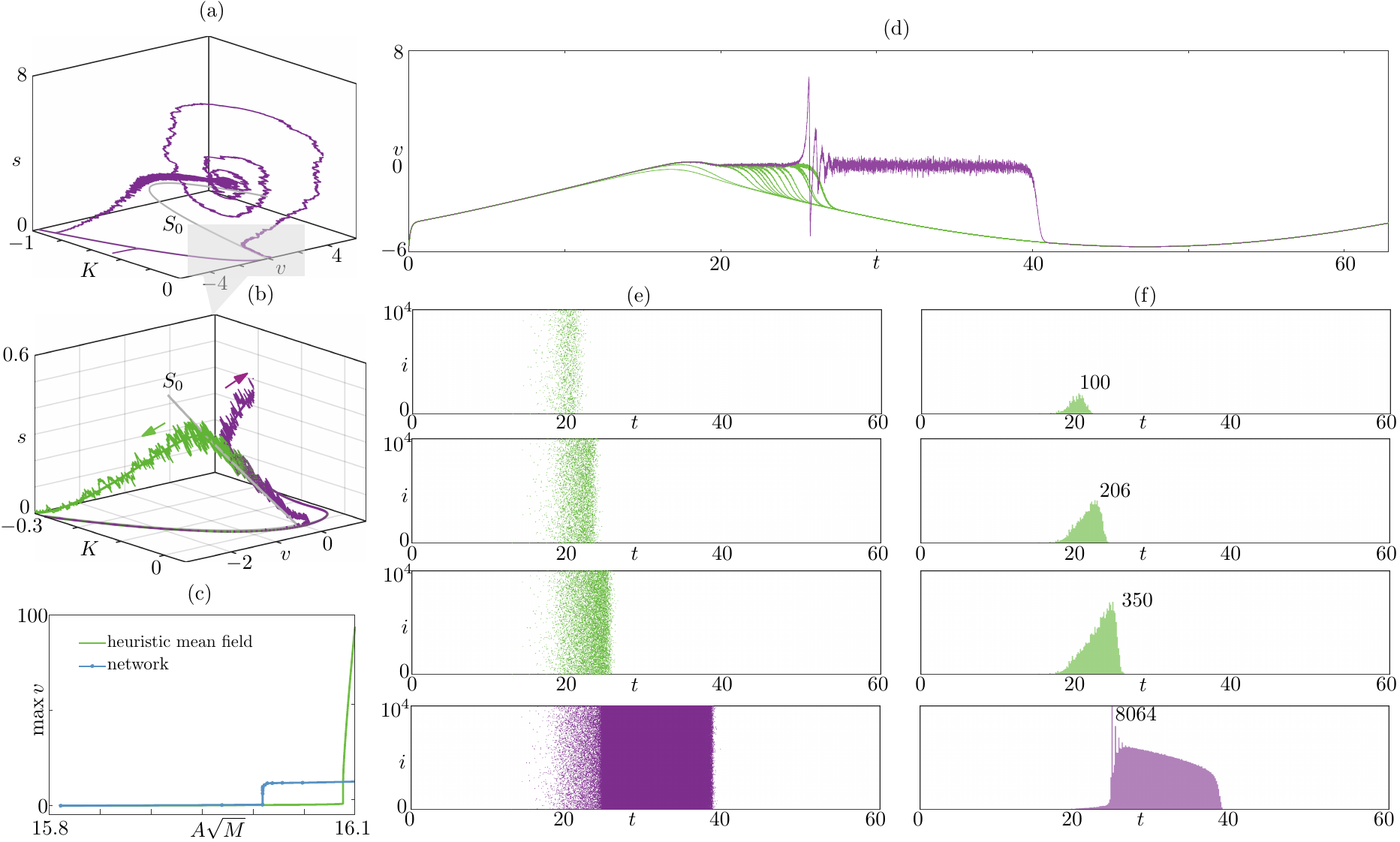}
  \caption{Dynamics of the sparse network~\eqref{eq:SparseNetw} with
    randomly-distributed currents and randomly-distributed connectivity. Panel (a) shows a network bursting solution plotted in the 3D phase space $(K,v,r)$ together
with the critical manifold $S_0$ of the heuristic mean-field
system~\eqref{eq:approxmeanField}. A down-up canard segment is visible in the
greyed out region. (b): A zoomed view of panel (a); in addition to $S_0$ and the bursting
solution of panel (a) in purple ($A=16.009453593274596/\sqrt{M}$), we show an orbit
with down-down canard segment obtained for slightly perturbed values of $A$, in
green ($A=16.009453593274599\sqrt{M}$); we also superimpose
solutions of the heuristic mean-field system, whose curves do not have fluctuations,
unlike the network ones. Panels (a,b) constitute numerical evidence that a down-to-up
canard-mediated transition to bursting exists in this network, as further shown in
panels (c--f). (c): Down-to-up route to bursting in the network (blue) and heuristic
mean field (green); besides some discrepancies discussed in the main text, both
branches contain a quasi-vertical segment typical of canard dynamics, bridging
between the non-bursting regime and the bursting one. (d) Representative network
solutions along the vertical branch in panel (d), 
between $A = 15.8113883008419/\sqrt{M}$ and $A =
16.0094535932746/\sqrt{M}$, 
displayed in the time series for
the voltage $v$; this solutions reveal a clear down-to-up transition, from non
bursting (green) to bursting (purple) orbits possessing canard segments (near $v=0$);
the peculiarity of these sparse-network canard solutions is that the associated rate
increases rapidly in the canard regime while the mean voltage remains approximately
constant, as evidenced in panels (e,f). (e): raster diagram of selected down-to-down
(green) and down-to-up (purple) network solutions from (d); the canard segments
manifests themselves in the raster plots: their onset coincides with the onset of
spiking, and their termination with the jump to the quiescent phase
(green) or the start of the tonic phase (purple); along the segment, the network
builds up rate (as shown in (f)). (f): histograms of firing events between $t$ and
$t+\Delta t$, with $\Delta t = 0.15$; along the canard segments the solution
increases the firing rate (the longer the segment, the higher the maximum rate in the
green histograms); the purple diagram is at a different scale with respect to the
green ones, and it represents a bursting solution. Parameters:
$N = 10^4$, $M=10^3$, $J =1$, $\tau_s = 0.015$, $\Delta_\gamma = 0.3$,
$\eps = 0.1$, $ \bar \eta = -0.5$, $\Delta = 10^{-4}$, $v_t = -v_r = 100$.
}
  \label{fig:sparse}
\end{figure*}

The slow-fast scenarios uncovered in the previous section are valid in a large variety
of (all-to-all coupled) QIF networks with exact mean-field limits. We now present
evidence that the phenomenon persists in sparse networks, for which no exact mean
field limit has been derived to date.

We present this extension for two main reasons: on one hand, we show that
the mechanism discussed in the previous section extends further, to sparse networks;
on the other, we want to emphasise that the availability of a mean-field description
is not strictly necessary for the canard phenomenon, which is supported by generic
network systems of QIFs with finite size. In \cref{sec:networkExc} we studied
networks with exact mean fields: since an ODE description was available for the case
$N\to \infty$, we used these ODEs to predict the region in parameter space where canard
dynamics occur, and to classify the folded singularities organising the transition
from non bursting to bursting patterns. However,
networks with finite size also support canard-mediated transitions, as evidenced in
\cref{fig:routeBurstN1} and \cref{fig:meanFieldFigs}, where the orbits are computed
for a very small ($N=1$) and a very large ($N=10^5$) network, respectively. 

Having a mean field description at our disposal is useful to pinpoint regions of
parameter space where canards will occur, through the study of $S_0$ and its folded
lines; also, large networks of neurons possess canard solutions that are almost
indistinguishable from their mean field ones. However, canard mediated transitions
are present (and can be documented) in finite-size networks, even when the mean field is
inexact, or unavailable in closed form.

To substantiate this claim, we study a sparse network of $N$ synaptically-coupled QIF
neurons. For a similar network, a heuristic mean field description has been proposed, based
on sparsity scaling arguments~\cite{DiVolo2018,bi2020}. The heuristic mean field is
in the form \eqref{eq:meanFieldGen} and hence canards of folded saddle type are
supported by this set of ODEs. However, this mean field is not the exact limit of a network of QIF
neurons, and the extent to which the mean field approximates the finite-size network
is also immaterial to find canards: as we shall see, both the heuristic mean field and the finite-size
network have canard-mediated routes to bursting, even if the two models do not agree
well in certain regions of parameter space.

We consider $N$ synaptically-coupled QIF neurons of the following form
\begin{equation}\label{eq:SparseNetw}
  \begin{aligned}
  &  V_i' = V_i^2 + \eta_i +  I(t) + \frac{J}{N\sqrt{M}}  \sum_{j=1}^N W_{ij} s_j,   \\
  & \tau_s s_i' =-s_i,
  \end{aligned}
\end{equation}
for $1\leqslant i\leqslant N$, where $M$ is an integer controlling the expected
number of connections of a neuron. More precisely, the connectivity matrix
 $W$, with entries $W_{ij}$, is a binary sparse matrix: the $i$th neuron receives
 input from $\gamma_i$ randomly selected neurons; the degree $\gamma_i$ is also
 random, given by
 \[
   \gamma_i =  \left\lfloor k_i \right\rfloor \chi_{[0,2M]}(k_i), 
   \qquad 
   k_i \stackrel{\rm i. i. d}{\sim} \frac{\Delta_\gamma \sqrt{M}}{(k-M)^2 +
   \Delta_\gamma^2 M},
 \]
 where $\chi$ is the indicator function. In practice, the connectivity matrix is
 established as follows: a candidate degree $k_i$ is extracted from a
 Cauchy distribution with
 center $M$ and HWHM $\Delta_\gamma \sqrt{M}$ and, if it lies in the interval $[0,2M]$, is
 rounded to the nearest lower integer to give $\gamma_i$; the $i$th row of the matrix
 $W$ has $\gamma_i$ randomly selected entries equal to 1, and the remaining
 $N-\gamma_i$ entries equal to $0$.

 In this model, the synaptic input scales as
 $1/\sqrt{M}$, the external forcing is given by $I(t) = A\sqrt{M}\sin(\epsi
 t)$, and the background currents $\eta_i$ are i.i.d, Cauchy distributed with peak at
 $\bar \eta \sqrt{M} $ and HWHM $\Delta \sqrt{M}$. In passing,
 we note that $\Delta_\gamma \neq \Delta$. Finally, the variables $v_i$ and $s_i$ are
 reset as in the other network examples presented above.

 With these scalings for the variable $M$ a heuristic, approximate mean-field
 description was
 proposed recently~\cite{DiVolo2018,bi2020}, for a system of inhibitory neurons with
 no forcing ($I(t) \equiv 0$), homogeneous currents ($\eta_i \equiv \eta$) and similar
 connectivity pattern. Reasoning in a similar fashion, we
 arrive at the following candidate approximate mean field given by
\begin{equation}\label{eq:approxmeanField}
 \begin{aligned}
 \epsi\dot{r} &= \Delta/\pi+2rv+J\Delta_\gamma s/\pi\\
 \epsi\dot{v} &= v^2+\sqrt{M}(K+Js)-(\pi r)^2\\
 \epsi\dot{s}  &= (-s+r)/\tau_s\\
 \dot{K} &=Q\\
 \dot{Q} &=-(K-\bar{\eta}).
 \end{aligned}
\end{equation}

Differently from~\cite{DiVolo2018,bi2020}, the model above has
heterogeneous currents, in addition to sparse, heterogeneous connectivity. Also we
consider an excitatory neuronal population, as opposed to a inhibitory one. The
inhibitory population considered in~\cite{DiVolo2018,bi2020}, with $J=1$ and a term
$-Js$ in the $v$-equation, is not suitable for studying excitability and transition
to bursting, as the critical manifold is not folded.

A bursting orbit for a network with $N=10^4$ neurons is visible in
\Cref{fig:sparse}(a), in the $(v,K,r)$-space. This simulation is done for a network
with sharply peaked current distribution ($\Delta \ll 1$), which generates a
particular bursting pattern, as we will now discuss.

Assuming that the heuristic mean-field description approximates the network
simulation, one can reason as in a standard MPR network: the burst is due to a family
of foci on the upper branch of an cubic-like critical manifold (visible in grey in
the figure). The figure shows that the manifold $S_0$ of the candidate mean field
captures well the geometry of the bursting orbit, and in particular it displays a
canard segment along the repelling branch of $S_0$.

A further inspection of $S_0$ reveals that, because $\Delta$ is small, both folds
$F^\pm$ of $S_0$ occur at vanishingly small values of $r$. We observe down-to-down and
down-to-up orbits in the network as well as in the heuristic mean field (see
\cref{fig:sparse}(b)), and this suggests the possibility of a continuous down-to-up
route to bursting.

It is important to note that the presence of nearby down-down and down-up solutions,
on its own, suffices to get a hint of the presence of canards; in this case, we also
have an approximating heuristic mean-field description with a computable manifold
$S_0$, which clearly helps us finding values of parameters where the nearby down-down
and down-up solution exist; the considerations that follow, however, hold for the
finite-size network, even if we obliterate $S_0$ from the pictures
\cref{fig:sparse}(a--b), and from the discussion above.

To uncover a continuous route to bursting \textit{in the network}, we compute several
orbits of the system when $A$ varies between $A = 15.8113883008419/\sqrt{M}$ and $A =
16.0094535932746/\sqrt{M}$; while $A$ changes, we keep the connectivity matrix $W$
constant, that is, we extract it once and reuse
it thereafter. The results are given in \cref{fig:sparse}(c--f). In \cref{fig:sparse}(c) we
show the voltage profiles. As anticipated, the canard structure of this network is
peculiar, because the currents are almost homogeneously distributed: the canard segment
in these orbits stretches along $v=0$; during this transition, however, the rate
increases sharply, as seen in the raster plot \cref{fig:sparse}(e) and in the
histograms in \cref{fig:sparse}(f). 
This means that the excitability threshold for this network occurs for states at
constant voltage and progressively large rate, unlike in the cases presented before.

In \cref{fig:sparse}(d) we compare the routes to bursting for the network and the
candidate mean field. We find a good agreement in the subthreshold regime, as well
as the presence of canard transitions in both systems, marked by quasi-vertical
branch segments. 

We also notice two types of discrepancies. Firstly, the value of $A$ at the
quasi-vertical segment differs slightly in the $N=10^4$ and in the mean field; this
is to be expected, given that the mean field is only heuristic. Solution types along
the canard explosion, however, are very similar in the two systems (not shown).
Secondly, we notice a discrepancy in the maximum voltage $v$, which is the solution
measure chosen for the network, especially in the bursting regime. This type of
discrepancy is not unexpected: in the finite-size network, each neuron has a voltage
at most equal to the reset value, hence the maximum mean voltage is capped in the network; the
heuristic mean field, on the other hand, supports solutions with a very large maximum
voltage. 

As anticipated, the discrepancy in \cref{fig:sparse}(c) is relevant if one
wants to assess the accuracy of the mean field in the bursting regime, but is
immaterial for the canard-mediated route to bursting: the blue route in
\cref{fig:sparse}(c), and the data in \cref{fig:sparse}(d--f) show such route
independently of the existence of $S_0$, or of the accuracy of the heuristic mean
field.

\section{Conclusions}
\label{sec:Conclusions}

%\textbf{Conclusion.}
The geometry of excitability and transition to bursting behavior in single neurons
and allied systems is governed by canard solutions, which act as thresholds and
determine the response of the system to slow parametric changes. We have shown that
this structure carries over in the mean-field limit of large populations of
excitable cells, as well as in large finite systems.  

In these cases, the average voltage of the population plays the same role as the
voltage in the single cell, and a well defined rate emerges as a new macroscopic
variable. If a separation of time scales exists between external input and voltage at
the level of a single cell, such separation persists at network level, between the
input and the coupled mean voltage and rate.

Two main results have been discovered at network level: (i) a large class of networks
of QIF neurons subject to an external stimulus with Cauchy-distributed background
currents undergoes a continuous transition to bursting (either up-to-down, or
down-to-up) upon increasing the forcing amplitude; to the best of our knowledge, this
statement holds for any QIF network amenable to the Ott-Antonsen reduction currently
derived in the literature, that is, expressible as \cref{eq:oasgen}. (ii) The
canard-mediated route to bursting is present also in sparse networks, for which
there is no exact limit; obtaining an approximate mean-field limit is convenient to
compare network trajectories to low-dimensional manifolds, but it is not strictly
necessary, and in fact continuous routes to bursting are present in small networks
too.

Since QIF neurons are general representatives of type-I neurons, we expect that
similar properties will survive in networks of more realistic cells, up to their mean
field-limit, which need not be an ODE. Introduction of inhibitory networks could
provide a connection between this work and the concept of balanced networks where,
depending on the details of the connection strengths, the firing is either driven by
the mean input (analogous to our tonic behavior) or by the fluctuations (analogous to
the excitable case). There is currently no available theory for such cases, except
for continuous coarse-grained networks~\cite{Avitabile:2017,avitabile2020local}, but
a numerical exploration of the averaged voltages and synaptic variables may reveal an
underlying low-dimensional structure, similarly to what has been found in
spatially-extended neural field models. 

% \nocite{*}
\bibliographystyle{apsrev4-1} 
\bibliography{references}
\end{document}